\documentclass{article}
\usepackage{latexsym, amsfonts, amsmath}
\newcommand{\be}{\begin{equation}}
\newcommand{\ee}{\end{equation}}
\newcommand{\bea}{\begin{eqnarray}}
\newcommand{\eea}{\end{eqnarray}}
\newcommand{\bean}{\begin{eqnarray*}}
\newcommand{\eean}{\end{eqnarray*}}

\newcommand{\brray}{\begin{array}}
\newcommand{\erray}{\end{array}}
\newcommand{\ben}{\begin{equation}{nonumber}}
\newcommand{\een}{\end{equation}{nonumber}}

\newtheorem{dfn}{Definition}[section]
\newtheorem{thm}[dfn]{Theorem}
\newtheorem{lmma}[dfn]{Lemma}
\newtheorem{ppsn}[dfn]{Proposition}
\newtheorem{crlre}[dfn]{Corollary}
\newtheorem{xmpl}[dfn]{Example}
\newtheorem{rmrk}[dfn]{Remark}
\newcommand{\bdfn}{\begin{dfn}}
\newcommand{\bthm}{\begin{thm}}
\newcommand{\blmma}{\begin{lmma}}
\newcommand{\bppsn}{\begin{ppsn}}
\newcommand{\bcrlre}{\begin{crlre}}
\newcommand{\bxmpl}{\begin{xmpl}}
\newcommand{\brmrk}{\begin{rmrk}}
\newcommand{\edfn}{\end{dfn}}
\newcommand{\ethm}{\end{thm}}
\newcommand{\elmma}{\end{lmma}}
\newcommand{\eppsn}{\end{ppsn}}
\newcommand{\ecrlre}{\end{crlre}}
\newcommand{\exmpl}{\end{xmpl}}
\newcommand{\ermrk}{\end{rmrk}}

\newcommand{\IC}{\mathbb{C}}

\newcommand{\IR}{\mathbb{R}}

\newcommand{\cla}{{\cal A}}
\newcommand{\clb}{{\cal B}}
\newcommand{\clc}{{\cal C}}
\newcommand{\cld}{{\cal D}}
\newcommand{\cle}{{\cal E}}
\newcommand{\clf}{{\cal F}}

\newcommand{\clh}{{\cal H}}
\newcommand{\cli}{{\cal I}}

\newcommand{\cll}{{\cal L}}
\newcommand{\clm}{{\cal M}}
\newcommand{\cln}{{\cal N}}
\newcommand{\clo}{{\cal O}}
\newcommand{\clp}{{\cal P}}
\newcommand{\clq}{{\cal Q}}

\newcommand{\cls}{{\cal S}}

\def\a*{{\cal A}_{h,*}}
\def\B{{\cal B}(h)}
\def\B1{{\cal B}_1(h)}
\def\b{{\cal B}^{\rm s.a.}(h)}
\def\b1{{\cal B}^{\rm s.a.}_1(h)}

\newcommand{\ot}{\otimes}

\newcommand{\raro}{\rightarrow}

\def \qed {$\Box$}

\def\a*{{\cal A}_{h,*}}
\def\B{{\cal B}(h)}
\def\B1{{\cal B}_1(h)}
\def\b{{\cal B}^{\rm s.a.}(h)}
\def\b1{{\cal B}^{\rm s.a.}_1(h)}

\begin{document}
\begin{center}
{\Large{\bf A note on geometric characterization of quantum isometries of classical manifolds  }}\\ 
{\large { \bf Debashish Goswami\footnote{Partially supported by Swarnajayanti Fellowship from D.S.T. (Govt. of India) and also acknowledges
 the Fields Institute, Toronto for providing hospitality for a brief stay when a small part of this work was done}, \bf Soumalya Joardar \footnote{Acknowledges support from CSIR}}}\\ 
Indian Statistical Institute\\
203, B. T. Road, Kolkata 700108\\
Email: goswamid@isical.ac.in\\
\end{center}
\begin{abstract}
If a compact quantum group acts isometrically on a (possibly disconnected) compact smooth Riemannin manifold $M$ such that the action (say $\alpha$) commutes with the Laplacian, i.e. isometric 
 in the sense of \cite{rigidity} then it is known (\cite{rigidity}) that the `differential' of the action preserves Riemannian inner product on forms in the sense that $<< (d \ot {\rm id})(\alpha(f)), 
  (d \ot {\rm id})(\alpha(g))>>=\alpha(<<df, dg>>)$ for all smooth functions $f,g$, where $<< \cdot, \cdot >>$ denotes the Riemannian inner-product viewed as a $C^\infty(M)$-valued
   inner-product on the bimodule of one-forms. In this note, we prove a partial converse to this, under the additional assumption that $M$ is oriented and the action preserves the orientation 
    in a suitable sense. Using this, an alternative line of arguments is given for the main result of \cite{rigidity}.

\end{abstract}
 \section{Introduction}
 It is  a very important and interesting problem in the theory of quantum groups and noncommutative geometry to study  `quantum symmetries' of various classical and quantum structures. 
 Indeed, symmetries of physical systems (classical or quantum) were conventionally modelled by group actions, and after the advent of quantum groups, group symmetries were naturally generalized to  symmetries given by quantum group action. In this context, it is  natural to think of quantum automorphism or the full quantum symmetry groups of various mathematical and physical structures. The underlying basic
 principle of defining a quantum automorphism group of a given mathematical structure consists of two steps : first, to identify (if possible) the group of automorphisms of
 the structure as a universal object in a
   suitable category, and then, try to look for the universal object in a similar but bigger category by replacing groups by quantum groups of appropriate type. 
   The formulation and study of such quantum symmetries in terms of universal Hopf algebras date back to Manin \cite{manin_book}.
      In the analytic set-up 
    of compact quantum groups, it was considered by S. Wang who defined and studied quantum permutation groups of finite sets and quantum automorphism groups 
    of finite dimensional algebras, such questions were taken up by a number of mathematicians including Banica, Bichon (see, e.g. \cite{ban_1}, \cite{ban_2}, \cite{bichon}, \cite{wang}), 
    and more recently in the framework of Connes' noncommutative geometry (\cite{con}) by Goswami, Bhowmick,  Skalski, Banica and others who have extensively studied 
    the quantum group of isometries  (or quantum isometry group) defined in \cite{Goswami}  (see also \cite{qorient}, \cite{qiso_trans}, \cite{qdisc} etc.). In this context, 
    it is important to compute such quantum isometry groups for classical (compact)  Riemannian manifolds. 

In the classical case, i.e. smooth group-actions on a Riemannian manifold $M$, the action commutes with the Hodge Laplacian $-d^*d$ if and only if its differential 
 is an isometry between co-tangent spaces for the given Riemannian structure, i.e. preserves the $C^\infty(M)$ valued inner product on the bimodule of smooth one-forms. It is natural to 
  see whether this extends to the quantum case. This is the aim of the present article. Indeed, it is easy to see one-way: if a compact quantum group action commutes with the Laplacian
   (this is what is termed as `isometric action' in \cite{Goswami}),  then it is smooth and preserves the inner product. However, we have been able to prove the converse only 
    in a slightly restricted set-up, namely when the manifold is oriented and the action also preserves the orientation in a suitable sense.

   \section{Notations and preliminaries}
  We follow the notations and set-up of \cite{rigidity}, which we briefly recall here. 
  All the
Hilbert spaces are over $\mathbb{C}$ unless mentioned otherwise. If $V$ is a
vector space over real numbers we denote its complexification by
$V_{\mathbb{C}}$. For a vector space $V$, $V^{'}$ stands for its algebraic dual.
For a C* algebra $\clc$, $\clm(\clc)$ will denote its
multiplier algebra. $\oplus$ and $\ot$ will denote the algebraic direct sum and
algebraic tensor
product respectively. We shall denote the $C^{\ast}$
algebra of bounded operators on a Hilbert space $\clh$ by $\clb(\clh)$ and the
$C^{\ast}$ algebra
of compact operators on $\clh$ by $\clb_{0}(\clh)$. $Sp$ ($\overline{Sp}$)
stands for the
linear span (closed linear span). 
Also WOT and SOT for the weak operator
topology and the strong operator topology respectively. Let $\clc$ be an
algebra. Then
$\sigma_{ij}:\underbrace{\clc\ot\clc\ot...\ot\clc}
_{n-times}\raro\underbrace{\clc\ot\clc\ot...\ot\clc}_{n-times}$ is the flip map
between $i$ and $j$-th place and
$m_{ij}:\underbrace{\clc\ot\clc\ot...\ot\clc}_{n-times}
\raro\underbrace{\clc\ot\clc\ot...\ot\clc}_{(n-1)-times}$ is the map obtained by
multiplying $i$ and $j$-th entry. In case we have two copies of an algebra we
shall simply denote by $\sigma$ and $m$ for the flip and multiplication map
respectively. We shall need several types of topological tensor products in
this paper: $\hat{\ot},\bar{\ot}, \bar{\ot}_{in}$ (to be explained
in subsequent sections). Also for a Hilbert space $\clh$ and a $C^{\ast}$
algebra or a locally
convex $\ast$ algebra $\clc$ we
shall consider the trivial Hilbert (bi)module $\clh\bar{\ot}\clc$ with the
obvious right and left action of $\clc$ on $\clh\bar{\ot}\clc$ coming from
algebra multiplication of $\clc$ and obvious $\clc$ valued inner product. When
$\clh=\mathbb{C}^{N}$, the (bi)module is called the trivial $\clc$ (bi)module of
rank $N$. Usually, we use $<,>$ and $<<,>>$ for the scalar valued inner product
(of a Hilbert space) and the algebra-valued inner product (of a Hilbert
bimodule) respectively. For a Hopf algebra $H$, for any $\mathbb{C}$-linear map
$f:H\raro H\ot H$, we write $f(q)=q_{(1)}\ot q_{(2)}$ (Sweedler's notation). For
an algebra or module $\cla$ and a $\mathbb{C}$-linear map $\Gamma:\cla\raro
\cla\ot H$, we shall also use an analogue of Sweedler's notation and write
$\Gamma(a)=a_{(0)}\ot a_{(1)}$.\\

We begin by recalling from \cite{Takesaki} the tensor product of two $C^{\ast}$
algebras $\clc_{1}$ and $\clc_{2}$ and let us choose the minimal or spatial
tensor product between two $C^{\ast}$ algebras. The corresponding $C^{\ast}$
algebra will be denoted by $\clc_{1}\hat{\ot}\clc_{2}$ throughout this paper.
However we need to consider more general topological spaces and algebras, namely, locally convex $\ast$-algebras embedded in 
 $C^*$ algebras with topology given by 
 countable family seminorms coming from closed derivations, which is a special class of $C^*$-normed smooth algebras in the sense of Blackadar and Cuntz (\cite{smooth}).
 \bdfn
 A unital Fr\'echet $\ast$ algebra $\cla$ will be called a `nice' algebra if there is a $C^*$-norm $\| \cdot \|$ on $\cla$ and the 
  the underlying locally convex topology of $\cla$ 
  comes from the family of seminorms $\{ \| \cdot \|_{\underline{\alpha}} \}$ given below:\\
$\{||x||_{\underline{\alpha}}=||\delta_{\underline{\alpha}}(x)||\}$, where
$\underline{\alpha}=(i_{1},...,i_{k}):1\leq i_{j}\leq k,k\geq 1$ is a multi
index or $\underline{\alpha}=\phi(null \ index)$,
$\delta_{\underline{\alpha}}=\delta_{i_{1}}..\delta_{i_{k}}$,
$\delta_{\phi}={\rm id}$ and where each $\delta_i$ denotes a $\| \cdot \|$-closable $\ast$-derivation from $\cla$ to itself.
\edfn

Given two such `nice' algebras
$\cla(\subset\cla_{1})$ and $\clb(\subset\clb_{1})$, where $\cla_1, \clb_1$ denote respectively the $C^*$-completion of $\cla,\clb$ in the corresponding $C^*$-norms,
 we
choose the injective tensor product norm on $\cla\ot \clb$, i.e. view it as a dense subalgebra of $ \cla_{1}\hat{\ot}\clb_{1}$. $\cla\ot\clb$ has
natural (closable $\ast$) derivations of the forms $\tilde{\delta}=\delta \ot {\rm id}$ as well as $\tilde{\eta}={\rm id} \ot \eta$ 
 where $\delta,\eta$ are closable $\ast$-derivations on $\cla$ and $\clb$ respectively. Clearly, $\tilde{\delta}$ commute with $\tilde{\eta}$.
 We topologize $\cla\ot \clb$ by the family
of seminorms coming from such derivations, i.e. $\{||.||_{\underline{\alpha}\underline{\beta}}\}$  where
$\| \cdot \|$ is the injective $C^*$-norm and 
$$ ||X||_{\underline{\alpha}\underline{\beta}}=||\tilde{\delta}_{\underline{\alpha}}\tilde{\eta}_{\underline{\beta}}(X)||,$$ $\underline{\alpha}
=(i_1,\ldots, i_k)$, $\underline{\beta}=(j_1,\ldots, j_l)$ some multi-indices as before and $\delta_i, \eta_j$'s being 
closable $\ast$-derivations on $\cla$ and $\clb$ respectively.
 
 We refer to \cite{rigidity} and references therein (in particular \cite{smooth}) for discussion on such algebras.  
 We actually need algebras of the form $C^\infty(M) \hat{\ot} \cla$ where $M$ a smooth compact 
   manifold possibly with boundary and $\cla$ is a unital $C^*$ algebra. By the nuclearity of $C^\infty(M)$ as a locally convex space
    the above tensor product is isomorphic with $C^\infty(M, \cla)$ with the Fr\'echet topology coming from the smooth 
     vector fields of $M$. Moreover, it is proved in \cite{rigidity} that $C^\infty(M, \cla)$ is again a nice algebra
     and we also have Fr\'echet $\ast$-algebra isomorphism between $C^\infty(M, \cla) \hat{\ot} \clb$ and $C^\infty(M, \cla \hat{\ot} \clb)$,
     where $\clb$ is any other unital $C^*$ algebra.

We also need Hilbert bimodules over such algebras and their internal and external tensor products.  
Let $\cle_{1}$ and $\cle_{2}$ be two Hilbert bimodules over two
locally
convex $*$ algebras(nice) $\clc_{1}$ and $\clc_{2}$ respectively. We denote the
algebra
valued inner product for the Hilbert bimodules by $<<,>>$. When the bimodule is
a
Hilbert space, we denote the corresponding scalar valued inner product by
$<,>$. Then $\cle_{1}\ot \cle_{2}$ has an obvious $\clc_{1}\ot \clc_{2}$
bimodule
structure, given by $(a\ot b)(e_{1}\ot e_{2})(a^{'}\ot b^{'})=ae_{1}a^{'}\ot
be_{2}b^{'}$ for $a,a^{'}\in \clc_{1}, b,b^{'}\in \clc_{2}$ and $e_{1}\in
\cle_{1},
e_{2}\in \cle_{2}$. Also define $\clc_{1}\ot \clc_{2}$ valued inner
product by $<<e_{1}\ot e_{2},f_{1}\ot f_{2}>>=<<e_{1},f_{1}>>\ot
<<e_{2},f_{2}>>$ for $e_{1}, f_{1}\in \cle_{1}$ and $e_{2},f_{2}\in \cle_{2}$.
We denote the completed module by $\cle_{1}\bar{\ot}\cle_{2}$. In fact
$\cle_{1}\bar{\ot}\cle_{2}$ is an
$\clc_{1}\hat{\ot}\clc_{2}$ bimodule. This is called the exterior tensor product
of two
bimodules. In particular if one of the bimodule is a Hilbert space $\clh$
(bimodule over $\mathbb{C}$) and the other is a $C^{*}$ algebra $\clq$
(bimodule over itself), then the exterior tensor product gives the usual
Hilbert $\clq$ module $\clh\bar{\ot}\clq$. When
$\clh=\mathbb{C}^{N}$, we have a natural identification of an element
$T=((T_{ij}))\in M_{N}(\clq)$ with the right $\clq$ linear map of
$\mathbb{C}^{N}\ot \clq$ given by
\begin{displaymath}
 e_{i}\mapsto \sum e_{j}\ot T_{ji},
\end{displaymath}
where $\{e_{i}\}_{i=1,...,N}$ is a basis for $\mathbb{C}^{N}$. We can take
tensor products of maps between two Hilbert bimodules under. We shall need such
tensor product of maps, which are `isometric' in some sense. Let
$T_{i}:\cle_{i}\raro\clf_{i}$, $i=1,2$ be two $\mathbb{C}$-linear maps and
$\cle_{i}$, $\clf_{i}$ be Hilbert bimodules over $\clc_{i}$, $\cld_{i}$
($i=1,2$) respectively. Moreover, suppose that
$<<T_{i}(\xi_{i}),T_{i}(\eta_{i})>>=\alpha_{i}<<\xi_{i},\eta_{i}>>$,
$\xi_{i},\eta_{i}\in \cle_{i}$ where $\alpha_{i}:\clc_{i}\raro\cld_{i}$ are
$\ast$-homomorphisms. Then it is easy to show that the algebraic tensor product
$T:=T_{1}\ot_{alg}T_{2}$ also satisfies
$<<T(\xi),T(\eta)>>=(\alpha_{1}\ot\alpha_{2})(<<\xi,\eta>>)$ and hence extends
to a well defined continuous map from $\cle_{1}\bar{\ot}\cle_{2}$ to
$\clf_{1}\bar{\ot}\clf_{2}$ again to be denoted by $T_{1}\ot T_{2}$. \\
\indent Let $\clb$, $\clc$, $\cld$ be three locally convex $\ast$ algebras. Also
let $\cle_{1}$ be an $\clb-\clc$ Hilbert bimodule and $\cle_{2}$ be a $\clc-\cld$
Hilbert bimodule. Then $\cle_{1}\ot_{\clc}\cle_{2}$ is an $\clb-\cld$ bimodule
in the usual way. We can
define a $\cld$ valued inner product that will make $\cle_{1}\ot_{\clc}\cle_{2}$
a pre-Hilbert $\clb-\cld$ bimodule. For that take $\omega_{1},\omega_{2}\in
\cle_{1}$ and $\eta_{1},\eta_{2}\in \cle_{2}$ and define\\
$$<<\omega_{1}\ot \eta_{1},\omega_{2}\ot \eta_{2}>>:=
<<\eta_{1},<<\omega_{1},\omega_{2}>>\eta_{2}>>.$$
\indent Let ${\cal I}=
\{\xi\in {\cal E}_{1}\ot_{\clc} {\cal E}_{2}$ such that $<<\xi,\xi>>=0 \}$.
Then define ${\cal E}_{1}\ot_{in} {\cal E}_{2}= {\cal E}_{1}\ot_{\clc} {\cal
E}_{2}/\cli$. We note that this semi inner product is actually an inner
product, so that $\cli=\{0\}$ (see proposition 4.5 of \cite{Lance}).
The topological completion of ${\cal E}_{1}\ot_{in} {\cal E}_{2}$ is called the
interior tensor product and we shall denote it by
$\cle_{1}\bar{\ot}_{in}\cle_{2}$. We denote the projection map from ${\cal
E}_{1}\ot_{\clc} {\cal E}_{2}$ to ${\cal E}_{1}\ot_{in} {\cal E}_{2}$ by $\pi$.
We also make the convention of calling a Hilbert $\cla-\cla$ bimodule simply
Hilbert $\cla$ bimodule.\\

\section{Compact quantum groups, their representations and actions}

\bdfn
A compact quantum group (CQG for short) is a  unital $C^{\ast}$ algebra $\clq$ with a
coassociative coproduct 
(see \cite{Van}) $\Delta$ from $\clq$ to $\clq \hat{\ot} \clq$  
such that each of the linear spans of $\Delta(\clq)(\clq\ot 1)$ and that
of $\Delta(\clq)(1\ot \clq)$ is norm-dense in $\clq \hat{\ot} \clq$.
\edfn
From this condition, one can obtain a canonical dense unital $\ast$-subalgebra
$\clq_0$ of $\clq$ on which linear maps $\kappa$ and 
$\epsilon$ (called the antipode and the counit respectively) are defined making the above subalgebra a Hopf $\ast$ algebra. In fact, this is  the algebra generated by the `matrix coefficients' of
the (finite dimensional) irreducible non degenerate representations (to be
defined 
 shortly) of the CQG. The antipode is an anti-homomorphism and also satisfies $\kappa(a^*)=(\kappa^{-1}(a))^*$ for $a \in \clq_{0}$.

\indent Let $\clh$ be a Hilbert space. Consider the multiplier algebra
$\clm(\clb_{0}(\clh)\hat{\ot} \clq)$. This algebra has two natural embeddings
into $\clm(\clb_{0}(\clh)\hat{\ot} \clq\hat{\ot} \clq)$. The first
one is
obtained by extending the map $x\mapsto x\ot 1$. The second one is obtained by
composing this map with the flip on the last two factors. We will write $w^{12}$
and $w^{13}$ for the images of an element $w\in \clm(\clb_{0}(\clh)\hat{\ot}
\clq)$
by these two maps respectively. Note that if $\clh$ is finite dimensional then
$\clm(\clb_{0}(\clh)\hat{\ot} \clq)$ is isomorphic to $\clb(\clh)\ot \clq$
(we don't need any topological completion).\\
\bdfn
Let $(\clq,\Delta)$ be a CQG. A unitary representation of $\clq$ on a Hilbert
space $\clh$ is a unitary element $\widetilde{U}\in \clm(\clb_{0}(\clh)\hat{\ot} \clq)$ such that $({\rm id}\ot \Delta)\widetilde{U}=\widetilde{U}^{12}\widetilde{U}^{13}$.
\edfn

Now Let $\clc$ be a nice algebra and $\clq$ be a compact quantum group. Then we have the following
\bdfn
A $\mathbb{C}$ linear map $\alpha:\clc\raro \clc\hat{\ot}\clq$ is
said to be a topological action of $\clq$ on $\clc$ if\\
1. $\alpha$ is a continuous $\ast$ algebra homomorphism.\\
2. $(\alpha\ot {\rm id})\alpha=({\rm id}\ot \Delta)\alpha$ (co-associativity).\\
3. $Sp \ \alpha(\clc)(1\ot \clq)$ is dense in $\clc\hat{\ot}\clq$ in the
corresponding Fr\'echet topology.
\edfn
Given a topological action $\alpha$, proceeding along the lines of \cite{Soltan}, we can prove the existence of a maximal dense $\ast$ subalgebra $\clc_{0}$ of $\clc$ such that the action is algebraic over $\clc_{0}$ in the sense that $\alpha(\clc_{0})\subset\clc_{0}\ot\clq_{0}$ and ${\rm Sp} \ \alpha(\clc_{0})(1\ot\clq_{0})=\clc_{0}\ot\clq_{0}$. Note that if the Fr\'echet algebra is a $C^{\ast}$ algebra, then the
definition of a topological action coincides with the usual $C^{\ast}$ action
of a compact quantum group.
\bdfn
A topological action $\alpha$ is said to be faithful if the $\ast$-subalgebra
of $\clq$ generated by the elements of the form $(\omega\ot {\rm id})\alpha$, where
$\omega$ is a continuous linear functional on $\clc$, is dense in $\clq$.
\edfn

We now generalize the notion of unitary representation on Hilbert
spaces to
 on Hilbert bimodules over nice, unital topological
$\ast$-algebras. Let $\cle$ be a Hilbert $\clc-\cld$ bimodule over  topological
$\ast$-algebras $\clc$ and $\cld$ and let $\clq$ be a compact quantum group. If
we consider
$\clq$ as a bimodule over itself, then we can form the exterior tensor product
$\cle\bar{\ot}\clq$ which is a
$\clc\hat{\ot} \clq-\cld\hat{\ot}\clq$ bimodule. Also let
$\alpha_{\clc}:\clc\raro \clc\hat{\ot} \clq$ and $\alpha_{\cld}:\cld\raro
\cld\hat{\ot} \clq$ be topological actions of $\clc$ and $\cld$ on $\clq$ in the
sense discussed earlier. Then using $\alpha$ we can
give $\cle\bar{\ot}\clq$ a $\clc-\cld$ bimodule structure given by
$a.\eta. a^{'}= \alpha_{\clc}(a)\eta \alpha_{\cld}(a^{'})$, for $\eta\in
\cle\bar{\ot}\clq$ and $a\in\clc,a^{'}\in \cld$ (but without any $\cld$ valued
inner
product).
\bdfn
A $\mathbb{C}$-linear map $\Gamma:\cle\raro \cle\bar{\ot}\clq$
is
said to be an $\alpha_{\cld}$ equivariant unitary representation of $\clq$ on
$\cle$ if\\
1. $\Gamma(\xi d)=\Gamma(\xi)\alpha_{\cld}(d)$ and
$\Gamma(c\xi)=\alpha_{\clc}(c)\Gamma(\xi)$) for $c\in\clc,d\in\cld$.\\
2. $<<\Gamma(\xi),\Gamma(\xi^{'})>>=\alpha_{\cld}(<<\xi,\xi^{'}>>)$, for
$\xi,\xi^{'}\in \cle$.\\
3. $(\Gamma\ot {\rm id})\Gamma=({\rm id}\ot \Delta)\Gamma$ (co associativity)\\
4. $\overline{Sp} \ \Gamma(\cle)(1\ot \clq)=\cle\bar{\ot}\clq$ (non
degeneracy).\\
\edfn
In the definition note that condition (2) allows one to define
$(\Gamma\ot {\rm id})$. 
We recall some relevant results from \cite{rigidity}.
\blmma
\label{int_eq}
Let $\cle_{1}$ be a Hilbert $\clb-\clc$ bimodule and $\cle_{2}$ be a
Hilbert $\clc-\cld$ bimodule. $\alpha_{\clb}$,$\alpha_{\clc}$,$\alpha_{\cld}$
be topological actions on a compact quantum group $\clq$ of topological
$\ast$-algebras $\clb,\clc,\cld$ respectively.
$\Gamma_{1}:\cle_{1}\raro\cle_{1}\bar{\ot} \clq$ and
$\Gamma_{2}:\cle_{2}\raro\cle_{2}\bar{\ot} \clq$ be $\alpha_{\clc}$ and
$\alpha_{\cld}$ equivariant unitary representations as discussed earlier. Then
$$<<\Gamma_{2}(\eta),<<\Gamma_{1}(\omega),\Gamma_{1}(\omega^{'})>>\Gamma_{2}
(\eta^ { ' })>>=\alpha_{\cld}<<\eta,<<\omega,\omega^{'}>>\eta^{'}>>.$$
\elmma

\blmma
\label{rep_lift_bimod}
Let $\cle_{1},\cle_{2},\clb,\clc,\cld,\alpha_{\clb},\alpha_{\clc},\alpha_{\cld},
\Gamma_{1},\Gamma_ {2},\clq$ be as in Lemma \ref{int_eq}.  Then we have an
$\alpha_{\cld}$ equivariant representation $\Gamma$ of $\clq$ on the Hilbert
$\clb-\cld$ bimodule $\cle_{1}\bar{\ot}_{in}\cle_{2}$. 
 \elmma

  In particular when $\cle$ is the trivial $\clc$-bimodule of rank $N$,
 we have the following:
\blmma
\label{unitary}
Given an $\alpha$ equivariant representation $\Gamma$ of $\clq$ on
$\mathbb{C}^{N}\ot\clc$ such that $\Gamma(e_{i}\ot
1_{\clc})=\sum_{j=1}^{N}e_{j}\ot b_{ji}$, $b_{ij}\in \clc\hat{\ot}\clq$ for
all $i,j=1,...,N$, where $\{e_{i};i=1,...,N \}$ is an orthonormal basis of
$\mathbb{C}^{N}$, then $U=((b_{ij}))_{i,j=1,....,N}$ is a unitary element of
$M_{N}(\clc\hat{\ot}\clq)$.
\elmma

\section{Locally convex $\ast$ algebras and Hilbert bimodules coming from
classical geometry}
 
Let $M$ be a compact, Riemannian manifold. As in \cite{rigidity}, we
 denote the algebra of real (complex) valued smooth
 functions on $M$ by $C^{\infty}(M)_{\mathbb{R}}$ $(C^{\infty}(M))$. Clearly
 $C^{\infty}(M)$ is the complexification of $C^{\infty}(M)_{\mathbb{R}}$ . It is a nice algebra whose locally convex topology is given by 
 a complete set of vector fields $\{\delta_{1},...,\delta_{N}\}$. For details of this topology we refer the reader to Subsection 5.1 of \cite{rigidity}.
 With this locally convex topology in fact $C^{\infty}(M)$ is a nice nuclear algebra so that we can consider topological action of $\clq$ on $C^{\infty}(M)$.
 Also let
$\Lambda^{k}(C^{\infty}(M))$ be the space of smooth $k$ forms on the manifold
$M$. We equip $\Lambda^{1}(C^{\infty}(M))$ with the natural locally convex
topology induced by the locally convex topology of $C^{\infty}(M)$ given by a
family of seminorms $\{p_{(U,(x_{1},\ldots,x_{n}), K, \beta)}\}$, where $
(U,(x_{1},\ldots,x_{n}) )$  is  a  local  cordinate chart, $\beta=(\beta_1,\beta_2,\ldots, \beta_r)$ is a multi-index with $\beta_i \in \{ 1,2,\ldots, n\}$ as before, 
 $K$ is  a compact 
subset, and 
$p_{(U,(x_{1},\ldots,x_{n}),K, \beta)}(\omega):= \sup_{x\in K, 1\leq i \leq n}|\partial_\beta f_{i}(x)|,$ where
$f_{i}\in C^{\infty}(M)$ such that $\omega|_{U}=\sum_{i=1}^{n}f_{i}dx_{i}|_{U}.$
It is clear from the definition that the differential map
$d:C^{\infty}(M)\raro \Omega^{1}(C^{\infty}(M))$ is Fr\'echet continuous.\\

Now for a $C^{*}$ algebra $\clq$, $\Lambda^{k}(C^{\infty}(M))\bar{\ot}
\clq$ has a natural $C^{\infty}(M)\hat{\ot}\clq$ bimodule structure. The left
action is given by $$(\sum_{i} f_{i}\ot
q_{i})(\sum_{j}[\pi_{(k)}(\omega_{j})]\ot q_{j}^{'})=
(\sum_{i,j}[\pi_{(k)}(f_{i}\omega_{j})]\ot q_{i}q_{j}^{'} )$$ 
The right action is similarly given. The inner product is given by
$$<<\sum_{i}\omega_{i}\ot q_{i},\sum_{j}\omega_{j}^{'}\ot
q_{j}^{'}>>=\sum_{i,j}<<\omega_{i},\omega_{j}^{'}>>\ot q_{i}^{*}q_{j}^{'}.$$
Topology on $\Lambda^{k}(C^{\infty}(M))\bar{\ot} \clq$ is given by requiring
$\omega_{n}\raro \omega$ if and only if $<<\omega_{n}- \omega,\omega_{n}-
\omega>>\raro 0$ in $C^{\infty}(M)\hat{\ot}\clq$ or $C^{\infty}(M,\clq)$. \\

\section{Hodge $\star$ map}
\indent Now consider the case when $M$ is orientable and a globally non-vanishing $n$-form ($n$ being the dimension on $M$)  has
been chosen. We introduce the Hodge star operator, which is a pointwise isometry
$\ast = \ast_{x}:\Lambda^{k}T^{\ast}_{x}M \rightarrow \Lambda^{n-k}T^{\ast}_{x}M
$. Choose a positively oriented orthonormal basis $\{\theta^{1}, \theta^{2},...,
\theta^{n} \}$ of $T^{\ast}_{x}M$. Sincs $\ast$ is a linear transformation it is
enough to define $*$ on a basis element $\theta^{i_{1}}\wedge
\theta^{i_{2}}\wedge...\wedge \theta^{i_{k}}(i_{1}<i_{2}<...<i_{k})$ of
$\Lambda^{k}T^{\ast}_{x}M$. Note that \\
\begin{eqnarray*}
{\rm dvol}(x)&=& \sqrt{{\rm det} (<\theta^{i},\theta^{j}>)}\theta^{1}\wedge \theta^{2}\wedge ...\wedge \theta^{n}\\
&=& \theta^{1}\wedge \theta^{2}\wedge ...\wedge \theta^{n}
\end{eqnarray*} 
\bdfn
$*(\theta^{i_{1}}\wedge \theta^{i_{2}}\wedge...\wedge \theta^{i_{k}})=\theta^{j_{1}}\wedge \theta^{j_{2}}\wedge...\wedge \theta^{j_{n-k}}$ where $\theta^{i_{1}}\wedge \theta^{i_{2}}\wedge...\wedge \theta^{i_{k}}\wedge\theta^{j_{1}}..\wedge \theta^{j_{n-k}}= {\rm dvol}(x)$.\\
\edfn
Since we are using $\mathbb{C}$ as the scalar field, we would like to define $\bar{\omega}$ for a $k$ form $\omega$. In the set-up introduced just before the definition we have some scalars $c_{i_{1},...,i_{k}}$ such that $\omega(x)=\sum c_{i_{1},...,i_{k}}\theta^{i_{1}}\wedge \theta^{i_{2}}\wedge...\wedge \theta^{i_{k}} $. Then define $\bar{\omega}$ to be $\bar{\omega}(x)=\sum \bar{c}_{i_{1},...,i_{k}}\theta^{i_{1}}\wedge \theta^{i_{2}}\wedge...\wedge \theta^{i_{k}}.$
Then the equation $<<\omega,\eta>>=*(\bar{\omega} \wedge *\eta)$ defines an
inner product on the Hilbert module $\Lambda^{k}(C^{\infty}(M))$ for all
$k=1,...,n$ which is the same as the $C^{\infty}(M)$ valued inner product
defined earlier. Then the Hodge star operator is a unitary between two Hilbert
modules $\Lambda^{k}(C^{\infty}(M))$ and $\Lambda^{n-k}(C^{\infty}(M))$ i.e.
$<<*\omega,*\eta>>= <<\omega,\eta>>$. Also for $\omega,\eta\in
\Lambda^{k}(C^{\infty}(M))$, we have
$\ast\omega\wedge\eta=<<\bar{\omega},\eta>>{\rm dvol}$. For further details about the
Hodge star operator we refer the reader to \cite{Rosen}.

 Hence we have $$(\ast\ot {\rm id}):\Lambda^{k}(C^{\infty}(M))\ot \clq\raro
\Lambda^{n-k}(C^{\infty}(M))\ot\clq.$$ Since Hodge $\ast$ operator is an
isometry, $(\ast\ot {\rm id})$ is continuous with respect to the Hilbert module
structure of $ \dot{\Lambda}(C^{\infty}(M))\hat\ot\clq$. So we have 
\begin{displaymath}
 (*\ot
{\rm id}):\Lambda^{k}(C^{\infty}(M))\bar{\ot} \clq\raro
\Lambda^{n-k}(C^{\infty}(M))\bar{\ot} \clq.
\end{displaymath}
We derive a characterization for $(*\ot
{\rm id}):\Lambda^{k}(C^{\infty}(M))\bar{\ot}\clq\raro
\Lambda^{n-k}(C^{\infty}(M))\bar{\ot}\clq $
 for all $k=1,...,n$.\\
\blmma
\label{star_char}
Let $\xi\in \Lambda^{n-k}(C^{\infty}(M))\bar{\ot}\clq$ and $X\in
\Lambda^{k}(C^{\infty}(M)) \bar {\ot}\clq$. Then the following are equivalent:\\
 (i) For all $Y\in
\Lambda^{k}(C^{\infty}(M))\bar{\ot}\clq$, 
\begin{eqnarray} \label{09134}
   \xi\wedge Y=<<\bar{X},Y>>({\rm dvol}\ot 1_{\clq}) 
   \end{eqnarray}
 (ii) $\xi=(*\ot {\rm id})X.$
\elmma
{\it Proof:}\\
$(i)\Rightarrow (ii)$:\\
Let $m\in M$. Choose a coordinate neighborhood $(U,x_{1},x_{2},....,x_{n})$
around $x$ in $M$ such that $\{dx_{1}(m),...,dx_{n}(m)\}$ is an orthonormal
basis for $T^{\ast}_{m}(M)$ for all $m\in U$. Now for any $l\in\{1,...,n\}$,
let $\Sigma_{l}$ be the set consisting of $l$ tuples $(i_{1},...,i_{l})$
such that $i_{1}<i_{2}<...<i_{l}$ and $i_{j}\in\{1,...,n\}$ for $j=1,...,l$.
For $I=(i_{1},...,i_{l})\in\Sigma_{l}$, we write $dx_{I}(m)$ for
$dx_{i_{1}}\wedge...\wedge dx_{i_{l}}(m)$. Also for
$I(=(i_{1},...,i_{p}))\in\Sigma_{p}$, $J(=(j_{1},...,j_{q}))\in\Sigma_{q}$, we
write $(I,J)$ for $(i_{1},...,i_{p},j_{1},...,j_{q})$.\\
\indent Now fix $I\in\Sigma_{k}$.Then we have a unique
$I^{\prime}\in\Sigma_{n-k}$ such that
\begin{displaymath}
 (\ast(dx_{I}))(m)=\epsilon(I)dx_{I^{\prime}}(m),
\end{displaymath}
where $\epsilon(I)$ is the sign of the permutation $(I,I^{\prime})$. Given
$X\in \Lambda^{k}(C^{\infty}(M))\bar{\ot}\clq$, for $m\in M$, we have
$q_{I}(m)\in\clq$ such that 
\begin{displaymath}
 X(m)= \sum_{I\in\Sigma_{k}}dx_{I}(m)q_{I}(m).
\end{displaymath}
Also for $\xi\in\Lambda^{n-k}(C^{\infty}(M))\bar{\ot}\clq$, we have
$w_{J}(m)\in\clq$ such that 
\begin{displaymath}
 \xi(m)= \sum_{J\in\Sigma_{n-k}}dx_{J}(m)w_{J}(m).
\end{displaymath}
Hence
\begin{displaymath}
 ((\ast\ot {\rm
id})X)(m)=\sum_{I\in\Sigma_{k}}\epsilon(I)dx_{I^{\prime}}(m)q_{I}(m),
\end{displaymath}
where $I^{\prime}\in\Sigma_{n-k}$ is as mentioned before.\\
\indent Now we fix some $L\in\Sigma_{k}$ and choose
$Y\in\Lambda^{k}(C^{\infty}(M))\bar{\ot}\clq$ such that
$Y(m)=dx_{L}(m)1_{\clq}$. Hence 
\begin{displaymath}
 (\xi\wedge Y)(m)= \sum_{J\in\Sigma_{n-k}}dx_{J}\wedge dx_{L}w_{J}(m).
\end{displaymath}
But for a fixed $L\in\Sigma_{k}$, there is a unique $J^{\prime}\in
\Sigma_{n-k}$ such that 
\begin{displaymath}
 dx_{J^{\prime}}(m)\wedge dx_{L}(m)= \epsilon(L){\rm dvol}(m).
\end{displaymath}
Hence 
\begin{displaymath}
 (\xi\wedge Y)(m)=\epsilon(L)w_{J^{\prime}}(m) {\rm dvol}(m).
\end{displaymath}
On the other hand
\begin{eqnarray*}
 && <<\bar{X}, Y>>(m){\rm dvol}(m)\\ 
&=& \sum_{I\in \Sigma_{k}}<dx_{I}(m)q_{I}(m)^{\ast},
dx_{L}(m)1_{\clq}>{\rm dvol}(m)\\ 
&=& q_{L}(m){\rm dvol}(m).
\end{eqnarray*}
Hence by (\ref{09234}), we have $q_{L}(m)=\epsilon(L)w_{J^{\prime}}$. So varying $Y$,
we have
\begin{displaymath}
 \xi(m) = \sum_{L\in\Sigma_{k}}\epsilon(L)dx_{J}(m)q_{L}(m),
\end{displaymath}
i.e. $(\ast\ot {\rm id})X=\xi.$
\indent The other direction of the proof is trivial.\\ \qed.

\section{ Smooth and inner-product preserving action}
\bdfn
A topological action of $\clq$ on the
Fr\'echet algebra $C^{\infty}(M)$ is called the smooth action of $\clq$ on the
manifold $M$. 
\edfn
It has been shown in \cite{rigidity} that 
a smooth action $\alpha$ of $\clq$ on $M$ extends to a $C^{\ast}$ action on $C(M)$ which is denoted by $\alpha$ again.

Moreover,  set $d\alpha (df):= (d\ot {\rm id})\alpha(f)$ for all $f\in C^{\infty}(M)$. The following is proved in \cite{rigidity}:\\
\bthm
\label{def_dalpha}
(i) $d\alpha$ extends to a well defined continuous map from
$\Omega^{1}(C^{\infty}(M))$ to $\Omega^{1}(C^{\infty}(M))\bar{\ot} {\clq}$
satisfying $d\alpha(df)= (d\ot {\rm id})\alpha(f)$. 
(ii) 
For every $x \in M$, the unital $\ast$-algebra $\clq_x$ generated by 
$(\nu\ot id)\alpha(f)(x)$, $\alpha(g)(x)$ with 
 $f,g\in C^{\infty}(M)$ and all smooth vector fields $\nu$ on $M$, is commutative.
\ethm

\bdfn

 We call a smooth action $\alpha$ on a Riemannian manifold $M$ to be inner-product preserving if  
\begin{eqnarray}
\label{inner_prod_pres_111}
 <<(d\ot {\rm id})\alpha (f),(d\ot {\rm id})\alpha (g)>>=\alpha(<<df,dg>> )
 \end{eqnarray}
 for all $f,g\in C^\infty(M)$. 
\edfn

In \cite{rigidity}, it is proved that an inner product preserving action induces a canonical unitary equivariant representation on each of the bimodules of forms
 $\Lambda^k(C^\infty(M))$, to be denoted by $d\alpha_{(k)}$, say. 
Moreover, using this, smooth actions on the total spaces of certain bundles $E^k_\epsilon$ have been constructed in Section 8 of \cite{rigidity}.   
 Let $T^*_\epsilon (M)$ be the total space of the cotangent bundle of a compact Riemannian manifold $M$ consisting of cotangent vectors of length less than or equal to $\epsilon$ for some positive epsilon. The arguments of Section 8 of 
 \cite{rigidity} go through verabtim to give a smooth action say $\eta$ on $C^\infty(T^*_\epsilon(M))$.
 
$T^{\ast}_{\epsilon}(M)$ is a compact $2n$ dimensional manifold. Note that
$\pi^{-1}(U)\cong
 U\times K$, where $K$ is an $n$-dimensional closed ball of radius $\epsilon$.
 Moreover $T^{\ast}_{\epsilon}(M)$ is orientable with the following natural
 orientation. At the point $(m,\omega)\in\pi^{-1}(U)$ and any choice
 $\omega_{1},...,\omega_{n}$ as before,
 ${\rm dvol}(m,\omega)\in\Lambda^{2n}(C^{\infty}(T^{\ast}_{\epsilon}(M)))$ is given
 by $(\omega_{1}\wedge\omega_{2}\wedge...\wedge\omega_{n}\wedge
 dt_{1}\wedge...\wedge dt_{n})(m,\omega)$. It can be seen to be independent of
 choice of $\omega_{1},...,\omega_{n}$ and also it is non zero everywhere.
 Henceforth, we shall consider $T^{\ast}_{\epsilon}(M)$ oriented with the
 globally defined non vanishing ${\rm dvol}$ as the choice of orientation.\\
\blmma
\label{orient}
The lifted action $\eta$  is
also
orientation preserving in the sense $d\eta_{(2n)}({\rm dvol})={\rm dvol}\ot 1_{\clq}$.
\elmma
{\it Proof}:\\
For $m\in M$, choose a trivializing neighborhood around $m$ and one forms
$\omega_{1},...,\omega_{n}$ such that $\{\omega_{1}(x),...,\omega_{n}(x)\}$
forms an orthonormal basis for $T^{\ast}_{x}(M)$ for all $x\in U$. Then there are $\clq$-valued functions
$f_{ij}$'s for $1\leq i,j \leq n$ such that $f_{ij}(m)\in\clq_{m}$ for all
$m\in M$ and $d\alpha(\omega_{i})(m)=\sum_{j}f_{ij}(m)\omega_{j}(m)$. Choose and fix some smooth non-negative function $\chi$ supported in $U$. By the commutativity of $\clq_m$, we get
\begin{eqnarray*}
 && d\alpha_{(n)}(\omega_{1}\wedge...\wedge\omega_{n})(m)\\
&=&\wedge_{j=1}^{n}(\sum_{j}f_{ij}\omega_{j})(m)\\
&=& \sum_{\sigma\in S_{n}}(sgn \
\sigma)f_{1\sigma(1)}(m)f_{2\sigma(2)}(m)...f_{n\sigma(n)}(m)(\omega_{1}
\wedge...\wedge \omega_{n})(m)\\
&=& \Delta(m)(\omega_{1}
\wedge...\wedge \omega_{n})(m).
\end{eqnarray*}
 where $\Delta(m)=det \left( \left(f_{ij}(m)\right)\right)$.
Also, we have 
\begin{eqnarray*}
 && d\eta_{n}(dt_{1}\wedge...\wedge dt_{n})(m,\omega)\\
&=& \wedge(\sum_{j}f_{ij}(m)dt_{j}(m,\omega))\\
&=& \Delta(m) (dt_{1}
\wedge...\wedge dt_{n})(m,\omega)
\end{eqnarray*}
Hence $d\eta_{(2n)}(\omega_{1}\wedge...\wedge \omega_{n}\wedge
dt_{1}\wedge...\wedge dt_{n})(m,\omega)=\Delta(m)^2
(\omega_{1}\wedge...\wedge \omega_{n}\wedge
dt_{1}\wedge...\wedge dt_{n})(m,\omega)$.\\
Now note that 
\begin{eqnarray*}
 &&\alpha(\chi)(m)^2<<d\alpha(\omega_{i}),d\alpha(\omega_{j})>>(m)\\
&=&\delta_{ij} \alpha(\chi)^2(m)
\end{eqnarray*}
 as $w_i$'s are orthonormal on the support of $\chi$.
 Moreover, each $f_{ij}(m)$ is self adjoint. Choosing any $\ast$-character
$\gamma$ on the commutative $C^{\ast}$ algebra $\clq_{m}$, we see that either $\gamma(\alpha(\chi)(m))=0$ or 
$((\gamma(f_{ij}(m))))$ is in $O_{n}(\mathbb{R})$ and its determinant $\gamma(\Delta(m)) $ is $ 1$ or
$-1$. Thus $\alpha(\chi)\Delta^2=\alpha(\chi)$, which implies 
\begin{displaymath}
 d\eta_{(2n)}(\chi {\rm dvol})=\alpha(\chi) ({\rm dvol} \ot 1),
\end{displaymath}
 and hence by a partition of 
 unity argument we complete the proof that $\eta$ is orientation-preserving.\qed\\

\section{Action commuting with the Laplacian , i.e. isometric}
\subsection{Isometric actions}

Recall the definition of $QISO^{\cll}$ for a spectral triple
 satisfying certain regularity conditions from \cite{rigidity}. In particular, all classical
spectral
 triples, i.e. those coming from the Dirac operator on the spinor bundle of a
 compact Riemannian spin manifold, do satisfy such conditions and hence
 $QISO^{\cll}$ is defined for them. In fact it easily follows from
\cite{Goswami}
 that
 one can go beyond spin manifolds and define (and prove existence of) such a
 quantum isometry group for any compact Riemannian manifold $M$ (without
 boundary)
 as the universal object in the category of CQG $\clq$ with a faithful action
 $\alpha$ on $C(M)$ such that $({\rm id}\ot \phi)\alpha(C^{\infty}(M))\subset
 C^{\infty}(M)$ for all state $\phi$ and commutes with the Hodge Laplacian (to
 be called the $L_{2}$ Laplacian) $\cll_{2}=-d^{*}d$ restricted to
 $L^{2}(M,{\rm dvol})$. We shall denote the universal object in this category by
 $QISO^{\cll}(M)$ in this paper. It is proved in (Theorem 3.8 of
\cite{qorient})
 that
 $QISO^{\cll}(M)\cong QISO^{+}_{I}(d+d^{*})$ where now $d$ is viewed as a map
on
 the Hilbert space of forms of all orders, i.e. the $L^{2}$ closure of
 $\oplus_{k=0}^{dim \ M}\Lambda^{k}(M)$.\\
 Furthermore it follows from the Sobolev theorem that $({\rm id}\ot
 \phi)\alpha(C^{\infty}(M))\subset C^{\infty}(M)$ for all state $\phi$. We have
the following (see \cite{rigidity} for a proof):\\
 \bthm
 \label{iso_smooth}
 $QISO^{\cll}$ (and hence any subobject in the category $\clq^{\cll}$) has a
 smooth action on $C^{\infty}(M)$. 
 \ethm

 Let us denote by $\cll$ the restriction of $\cll_{2}$ to
 $C^{\infty}(M)$, viewed as a Fr\'echet continuous operator (to be called the
 `geometric Laplacian'). When $M$ is oriented we can also write it as
 $-(*d)^{2}$, where $*$ is the Hodge * operator as discussed in subsection 5.2.
 As $C^{\infty}(M)$ is a core for $\cll_{2}$, it is clear that a CQG action
 $\alpha:C(M)\raro C(M)\hat{\ot} \clq$ is isometric (i.e. $(\cla,\alpha)$ is an
 object in $\clq^{\cll}$) if and only if $\alpha$ is smooth and commutes with
 $\cll$ in the sense that $\alpha\circ\cll=(\cll\ot 1)\alpha$.\\
 For the purpose of this paper, we need to extend the above formulation of
 quantum isometry group to manifolds with boundary. Choosing the Dirichlet
boundary
 condition, we take $d$ to be the closure of the unbounded operator with domain
 $\clc=\{ f\in C^{\infty}(M):f|_{\partial M=0} \}$. 
 \bdfn
 For a
 compact manifold with boundary we call a smooth action
 $\alpha:C(M)\raro C(M)\hat{\ot}\clq$ to be isometric if it maps $\clc$ into
 $\clc
 \hat{\ot} \clq$ and commutes with $\cll_{2}$ on $C^{\infty}(M)$.
 \edfn
 \brmrk 
 For a manifold with boundary, commutation with the geometric Laplacian $\cll$
 may not be sufficient to imply that $\alpha$ is isometric. We also require
the
 condition that $\alpha(\clc)\subset \clc \hat{\ot}\clq$. We can prove the
existence of $QISO^{\cll}$ as well as the smoothness of the action of
$QISO^{\cll}$ as in \cite{Goswami}. It is a consequence of the
 fact that the Dirichlet Laplacian has discrete spectrum with finite
dimensional
 eigen spaces and the estimate $||e_{j}(f)||_{\infty}\leq
C\lambda_{j}^{\frac{n-1}{2}}||f||_{2}$ of the eigen vectors of the Laplacian
(see page
9 of \cite{xu}).
 \ermrk

 \blmma
 \label{inner_prod_pres}
 If $\alpha$ commutes with the geometric Laplacian $\cll$ on $\cla$, then 
  $\alpha$ is inner product preserving.
 \elmma 
 {\it Proof:}\\
 \begin{eqnarray*}
  && <<(d\ot {\rm id})\alpha (f),(d\ot {\rm id})\alpha (g)>>\\
  &=& <<df_{(0)},dg_{(0)}>>\ot f_{(1)}^{\ast}g_{(1)}\\
  &=&
 [\cll(\overline{f_{(0)}}g_{(0)})-\cll(\overline{f_{(0)}})g_{(0)}-\overline{f_{
 (0) } } \cll(g_{(0)})]\ot f_{(1)}^{\ast}g_{(1)}\\
 \end{eqnarray*}
 On the other hand
 \begin{eqnarray*}
  && \alpha(<<df,dg>>)\\
  &=& \alpha [\cll(\bar{f}g)-\cll(\bar{f})g-\bar{f}\cll(g)]\\
  &=&
 [\cll(\overline{f_{(0)}}g_{(0)})-\cll(\overline{f_{(0)}})g_{(0)}-\overline{f_{
 (0) } }
 \cll(g_{(0)})]\ot f_{(1)}^{\ast}g_{(1)}( \ since \ \alpha \ commutes \ with
 \ \cll)\\
 \end{eqnarray*}
 \qed\\
\subsection{Geometric characterization of orientation-preserving isometric action}
Our aim of this subsection is to prove  a partial converse to the fact that an isometric action is inner product preserving. More precisely, we shall prove the following
\blmma 
\label{main_lemm_comm}
Let $N$ be an $m$-dimensional compact, oriented, Riemannian manifold (possibly
with boundary) with ${\rm dvol}\in \Lambda^{m}(C^{\infty}(N))$ a globally defined
nonzero form. Moreover let $\eta$ be a smooth inner product preserving action on
$N$  such that $d\eta_{(m)}({\rm dvol})={\rm dvol}\ot 1$. Then $\eta$ commutes with the
geometric Laplacian.
\elmma
{\it Proof}:\\
First we note that as $\eta$ is an inner product preserving smooth action, by
the results of \cite{rigidity} (Corollary 7.12) it lifts to an $\alpha$-equivariant unitary
representations
$d\eta_{(k)}:\Lambda^{k}(C^{\infty}(N))\raro
\Lambda^{k}(C^{\infty}(N))\bar{\ot}\clq$ for all $k=1,...,m$. Note that without
loss of generality we can replace ${\rm dvol}$ by
$\frac{{\rm dvol}}{<<{\rm dvol},{\rm dvol}>>^{\frac{1}{2}}}$ and assume that $<<{\rm dvol},{\rm dvol}>>=1$,
since if $d\eta_{(m)}$ preseves ${\rm dvol}$, it also preserves the normalized ${\rm dvol}$.
First we claim
that \be \label{09234} 
\forall \ k=1,...,m,  d\eta_{(m-k)}(*\omega)\wedge \beta=
 <<\overline{d\eta_{(k)}(\omega)},\beta>>({\rm dvol}\ot 1_{\clq})\ee $ \forall \omega \in \Lambda^k(C^\infty(N))$, 
$ \forall \beta \in
 \Lambda^{k}(C^{\infty}(N))\bar{\ot}\clq.$
For that let $\beta=d\eta_{(k)}(\omega^{'})(1\ot q^{'})$. Then 
\begin{eqnarray*}
&& d\eta_{(m-k)}(*\omega)\wedge \beta \\
&=&d\eta_{(m-k)}(*\omega)\wedge d\eta_{(k)}(\omega^{'})(1\ot q^{'})\\
&=& d\eta_{(m)}(\ast(\omega)\wedge \omega^{\prime})(1\ot q^{\prime})\\
&=& d\eta_{(m)}(<<\bar{\omega},\omega^{\prime}>>{\rm dvol})(1\ot q^{\prime})\\
&=&\eta<<\bar{\omega},\omega^{'}>> ({\rm dvol}\ot q^{'}) \ 
\end{eqnarray*}
On the other hand from unitarity of $d\eta_{(k)}$,\\
\begin{eqnarray*}
&& <<\overline{d\eta_{(k)}(\omega)},d\eta_{(k)}(\omega^{'})(1\ot q^{'})>>\\
&=&\eta(<<\bar{\omega},\omega^{'}>>)(1\ot q^{'}).
\end{eqnarray*}
So by replacing $\beta$ by finite sums of the type
$\sum_{i}d\eta_{(k)}(\omega_{i})(1\ot q_{i})$, we can show that for $\omega \in
\Lambda^{k}(C^{\infty}(N))$ and $\beta\in  Sp \
d\eta_{(k)}\Lambda^{k}(C^{\infty}(N))(1\ot \clq)$,
$$d\eta_{(m-k)}(*\omega)\wedge \beta=<<d\eta_{(k)}(\bar{\omega}),\beta>>({\rm dvol}\ot
1_{\clq}).$$
Now, since Sp $d\eta_{(k)}(\Lambda^{k}(C^{\infty}(N)))(1\ot
\clq)$ is dense in
$\Lambda^{k}(C^{\infty}(N))\bar{\ot}\clq$, we get a
sequence $\beta_{n}$
belonging to Sp
$d\eta_{(k)}(\Lambda^{k}(C^{\infty}(N)))(1\ot \clq)$ such
that
$\beta_{n}\raro \beta$ in the Hilbert module
$\Lambda^{k}(C^{\infty}(N))\bar{\ot}\clq$.\\
But we have  
 $$d\eta_{(m-k)}(*\omega)\wedge
\beta_{n}=<<d\eta_{(k)}(\bar{\omega}),\beta_{n}>>({\rm dvol}\ot 1_{\clq}).$$
 Hence the claim follows from the continuity of $<<,>>$ and $\wedge$ in the
Hilbert module $\dot{\Lambda}(C^{\infty}(N)\bar{\ot}
\clq$.\\
 
  Combining Lemma \ref{star_char} and (\ref{09234}) we immediately
conclude the following:\\
 \begin{eqnarray}
 \label{17}
d\eta_{(m-k)}(*\omega)=(*\ot {\rm id})d\eta_{(k)}(\omega) \ for \ k\geq 0.
 \end{eqnarray}
Now we can prove that $\eta$ commutes with the geometric Laplacian of
$N$.
 For $\phi\in C^{\infty}(N)$,
   \begin{eqnarray*}
   &&\eta(\ast d\ast d \phi)\\
   &=& \ (\ast\ot {\rm id})d\eta_{(m)}(d\ast d\phi) \ ({\rm by} ~ {\rm
equation}~ (\ref{17})~{\rm with}~k=m)
 \\
   &=& \ (\ast d\ot {\rm id})d\eta_{(m-1)}(\ast d\phi)\\
   &=& \ (\ast d \ot {\rm id})(\ast\ot {\rm id})d\eta(d\phi) \ ({\rm again} ~{
\rm by }~{\rm  equation}~
(\ref{17}))
 \\
   &=&  (\ast d \ot {\rm id})(\ast d\ot {\rm id})\eta(\phi) \\
   &=&  ( (\ast d)^2 \ot {\rm id})\eta(\phi).\\
   \end{eqnarray*}
\qed\\


  \section{Application: Non existence of genuine CQG action}
\subsection{The stably parallelizable case}
   We now introduce the notion of stably parallelizable manifolds.
  \bdfn
  A manifold $M$ is said to be stably parallelizable if its tangent bundle is stably trivial.
  \edfn
  We recall the following from \cite{Singhof}:\\
  \bppsn
  A manifold $M$ is stably parallelizable if and only if it has trivial normal
  bundle when embedded in a Euclidean space of dimension higher than twice the
  dimension of $M$.
  \eppsn
  {\bf Proof}: see discussion following the Theorem (7.2) of \cite{Kol}.\\
  \qed\\
  We note that parallelizable manifolds (i.e. which has trivial tangent bundles)
  are in particular stably parallelizable. Moreover, given any compact Riemannian
  manifold $M$, its orthonormal frame bundle $O_{M}$ is parallelizable. Also given any stably parallelizable manifold $M$, the total space of its cotangent bundle is again stably parallelizable.\\
    Recall the manifold $T^{\ast}_{\epsilon}(M)$ for a smooth compact manifold $M$. If $M$ is isometrically embedded in some ${\mathbb R}^{N}$, then we consider
 the set $V=\{(u,v):u\in{\mathbb R}^{N}, v\in {\mathbb R}^{N} {\rm such \ that} \ ||v||\leq\epsilon\}$, then we define $\Phi:T^{\ast}_{\epsilon}(M)\raro V$ by 
 $\Phi(m,v)=(\phi(m),d\phi(v))$ where $\phi$ is the isometric embedding of $M$. Then it is easy to see that $\Phi(T^{\ast}_{\epsilon}(M))$ is a submanifold of
 $V$ and in fact it is a neat submanifold of $V$. So by Theorem 6.3 (page 114) of \cite{Morris}, we have 
 \blmma
 \label{tub_neat}
 $T^{\ast}_{\epsilon}(M)$ has a tubular neighborhood for some $\delta>0$ in $V$. It is denoted by $\cln_{\delta}(T^{\ast}_{\epsilon}(M))$.
\elmma
As $T^{\ast}_{\epsilon}(M)$ has a trivial normal bundle in $V$, the tubular neighborhood is actually diffeomorphic to $T^{\ast}_{\epsilon}(M)\times B_{\delta}^{2(N-n)}$, where $n$ is the dimension of the manifold. We denote the global coordinates for $B_{\delta}^{2(N-n)}$ by $u_{1},...,u_{2(N-n)}$. Let us recall from \cite{rigidity} the lift of a smooth action on a compact, stably parallelizable, Riemannian manifold to its tubular neighborhood. We denote the lift of the action $\eta$ on the manifold $T^{\ast}_{\epsilon}(M)$ to $\cln_{\delta}(T^{\ast}_{\epsilon}(M))$ by $\eta^\prime$. Also we take the canonical volume form of $\cln_{\delta}(T^{\ast}_{\epsilon}(M))$ to be ${\rm dvol}\ot du_{1}\wedge...\wedge du_{2(N-n)}$, where ${\rm dvol}$ is the volume form of $T^{\ast}_{\epsilon}(M)$. Then we have
\bppsn
\label{inn_orient}
(i) $\eta^\prime$ is inner product preserving.\\
(ii) $\eta^{\prime}$ is orientation preserving.
\eppsn
{\it Proof}:\\
The first statement was proved in \cite{rigidity} (Lemma 9.3). For the second statement it is enough to observe that $\eta^{\prime}(u_{i})=u_{i}$ for all $i=1,...,2(N-n)$ and for the functions of the form $f\circ \pi$, where $\pi$ is the projection of the normal bundle, $\eta^\prime(f\circ \pi)(y)=\eta(f)(\pi(y))$ for all $y\in \cln_{\delta}(T^{\ast}_{\epsilon}(M))$, which follows from the definition of the extension $\eta^\prime.$ and the fact that $\eta$ is orientation preserving. \qed\\
 Let $\{y_{i}:i=1,...,N\}$ be the
standard coordinates for $\mathbb{R}^{N}$.
We will also use the same notation for the restrictions of $y_{i}$'s if no
confusion arises.\\
  \bdfn
  A twice continuously differentiable, complex-valued function $\Psi$ defined on
a non empty, open set $\Omega\subset \mathbb{R}^{N}$  is said to be 
harmonic on $\Omega$ if $$\cll_{\mathbb{R}^{N}} \Psi\equiv 0,$$
  where $\cll_{\mathbb{R}^{N}}\equiv \sum_{i=1}^{N}\frac{\partial^{2}}{\partial y_{i}^{2}}$.\\
  \edfn
 
   \blmma
\label{quad}
Let $W$ be a manifold (possibly with boundary) embedded in some $\mathbb{R}^{N}$
and $\{y_{i}\}$'s for $i=1,...,N$, be the coordinate functions for
$\mathbb{R}^{N}$ restricted to $W$. If $W$ has non empty interior
in $\mathbb{R}^{N}$, then $\{1,y_{i}y_{j},y_{i}:1\leq i,j\leq N\}$ are
linearly independent, i.e. $\{ 1, y_1, \ldots, y_N\}$ are quadratically independent.
\elmma

\indent We call any action which preserves $V=\{1,y_{1},...,y_{N}\}$ affine.
  \blmma
  \label{affine}
 Let $\Phi$ be a smooth action of a CQG on a compact subset of $\mathbb{R}^{N}$
which commutes with $\cll_{\mathbb{R}^{N}}$, Then $\Phi$ is
affine
i.e.
\begin{displaymath}
 \Phi(y_{i})=1\ot q_{i}+ \sum_{j=1}^{N}y_{j}\ot q_{ij}, \ for \ some \
q_{ij},q_{i}\in \clq,
\end{displaymath}
 for all $i=1,...,N$, where $y_{i}'$s are coordinates of
$\mathbb{R}^{N}$.
  \elmma
 {\it Proof:}\\
 As $\Phi$ commutes with the geometric Laplacian and
$\cll_{\mathbb{R}^{N}}\frac{\partial}{\partial y_{j}}=\frac{\partial}{\partial
y_{j}}\cll_{\mathbb{R}^{N}}$, $\cll_{\mathbb{R}^{N}}y_{j}=0$ for all $j$,
 we get
 \begin{eqnarray*}
 && (\cll_{\mathbb{R}^{N}}\ot {\rm id})(\frac{\partial}{\partial
y_{j}}\ot id)\Phi(y_{i})\\
 &=& (\frac{\partial}{\partial y_{j}}\ot
{\rm id})\Phi(\cll_{\mathbb{R}^{N}} y_{i})\\
 &=& 0.
 \end{eqnarray*}
Let $D_{ij}(y)=((\frac{\partial}{\partial y_{i}}\ot {\rm id})\Phi(y_{j}))(y)$. Note
that as $d\Phi$ is an $\Phi$-equivariant unitary representation, by Lemma
\ref{unitary} $((D_{ij}(y)))_{i,j=1,...,N}$ is unitary for all $y\in
W$.
Pick $y_{0}$ in the interior of $W$(which is non empty). Then the new $\clq$ valued matrix $((G_{ij}(y)))=((D_{ij}(y)))((D_{ij}(y_{0})))^{-1}$ is unitary (since $D_{ij}(y)$ is so).\\
$G_{ij}(y)$ is unitary for all $y$ $\Rightarrow$ $|\psi(G_{ij}(y))|\leq 1$   
And $|\psi(G_{ii}(y_{0}))|=1$.
 $\psi(G_{ii}(y))$ is a harmonic function on an open {\it connected} set $Int
(W)$ which attains its supremum at an interior point. Hence by
corollary 1.9 of \cite{Harmonic} we conclude that
$\psi(G_{ii}(y))=\psi(G_{ii}(y_{0}))=1$. $((G_{ij}(y)))$ being unitary for all
$y$, we get $G_{ij}=\delta_{ij}.1_{\clq}$. Then
$((D_{ij}(y)))((D_{ij}(y_{0})))^{-1}=1_{M_{N}(\clq)}$,
i.e. $((D_{ij}(y)))=((D_{ij}(y_{0})))$ for all $y\in W$.
Hence $\Phi$ is affine with $q_{ij}=D_{ij}(y_{0})$\qed  \\
We also state the following Lemma without proof. For the proof reader might see \cite{rigidity}.
\blmma
\label{noaction}
Let $\clc$ be a unital commutative $C^{*}$ algebra and $x_1,x_2,\ldots, x_N$ be
self adjoint elements of $\clc$ such that $\{  x_{i}x_{j}:1\leq i\leq j\leq N\}$
are linearly independent and $\clc$ be a unital $C^{*}$ algebra generated by
$\{x_{1},x_{2},...,x_{N}\}$. Let $\clq$ be a compact quantum group
acting faithfully on $\clc$ such that the action  leaves the span of 
$\{x_{1},x_{2},...,x_{N}\}$ invariant. Then $\clq$ must be commutative
as a $C^{*}$ algebra, i.e. $\clq\cong C(G)$ for some compact group $G$.
\elmma

\brmrk
This is the only place where we need the manifold to be connected.
\ermrk
\qed\\
\bcrlre
\label{main_para}
Let $M$ be a smooth, compact, orientable, connected, stably parallelizable
manifold. Then if $\alpha$ is a faithful smooth action of a CQG $\clq$. 
 Then $\clq$ must be commutative as a
$C^{\ast}$ algebra i.e.
$\clq\cong C(G)$ for some compact group $G$.
\ecrlre
{\it Proof}:\\
 First recall from \cite{rigidity} (Theorem 7.13) that given a smooth action $\alpha$ of a CQG $\clq$ on a compact Riemannian manifold $M$, we can equip the manifold with a Riemannian structure such that the action becomes inner product preserving. So, by applying the averaging trick we reduce the action to an inner product preserving action first. Then we lift the action to the total space of the cotangent bundle.
 By Lemma \ref{orient}, the action is also orientation preserving. Now again using the averaging trick we equip the total space of the cotangent bundle with a new 
 Riemannian metric such that the action is inner product preserving. Then we lift this orientation preserving and inner product preserving action to the 
 tubular neighborhood of the total space of the cotangent bundle (which exists by Lemma \ref{tub_neat}) and by Proposition \ref{inn_orient}, we see that it is still orientation and inner product preserving. So by lemma \ref{main_lemm_comm}, it commutes with the geometric Laplacian of the tubular neighborhood, which is an open subset of ${\mathbb R}^{N}$ for some $N$. Now by applying Lemma \ref{affine}, Lemma \ref{quad} and Lemma
\ref{noaction}, we complete the proof.\qed\\
Using the above result and using the isometric lift of an action on a manifold to the total space of its orthonormal frame bundle (which is parallelizable) we get the main result of 
\cite{rigidity} which states that

\bthm
\label{main_res}
Let $\alpha$ be a smooth, faithful action of a CQG $\clq$ on a compact, connected
smooth manifold $M$. 
Then $\clq$ must be
commutative as a $C^{\ast}$ algebra i.e. $\clq\cong C(G)$ for some compact
group $G$.
\ethm
 \bcrlre
 The quantum isometry group of a compact, connected, Riemannian manifold
 coincides with the classical isometry group of the manifold.
\ecrlre
{\it Proof}:\\
Follows from the fact (Theorem \ref{iso_smooth}) that an isometric action of a
compact quantum group is smooth. 
\qed\\

\end{document}